\documentclass[12pt]{amsart}
\usepackage{latexsym,amsmath,amssymb,epsfig,amsthm}
\usepackage{graphicx}
\usepackage{tikz}
\usetikzlibrary{calc}
\usepackage[enableskew,vcentermath]{youngtab}
\usepackage{hyperref}
\hypersetup{colorlinks=false}
\input epsf
\newtheorem{theorem}{Theorem}[section]
\newtheorem{proposition}[theorem]{Proposition}

\newtheorem{conjecture}[theorem]{Conjecture}

\newtheorem{example}[theorem]{Example}
\newtheorem{lemma}[theorem]{Lemma}
\newtheorem{remark}[theorem]{Remark}

\newcommand{\asc}{{\rm asc}}

\newcommand{\des}{{\rm des}}
\newcommand{\Des}{{\rm Des}}
\newcommand{\Ehr}{{\rm Ehr}}

\newcommand{\rev}{{\rm rev}}

\newcommand{\cC}{{\mathcal C}}
\newcommand{\eE}{{\mathcal E}}

\newcommand{\hH}{{\mathcal H}}
\newcommand{\gG}{{\mathcal G}}

\newcommand{\lL}{{\mathcal L}}
\newcommand{\mM}{{\mathcal M}}

\newcommand{\qQ}{{\mathcal Q}}

\newcommand{\wW}{{\mathcal W}}
\newcommand{\zZ}{{\mathcal Z}}

\newcommand{\bu}{{\mathbf u}}

\newcommand{\RR}{{\mathbb R}}

\newcommand{\NN}{{\mathbb N}}
\newcommand{\ZZ}{{\mathbb Z}}

\renewcommand{\to}{\rightarrow}

\newcommand{\sm}{{\smallsetminus}}
\begin{document}
\title[Polytopes associated to integer compositions]
{Lattice point enumeration of polytopes associated 
to integer compositions}

\author{Christos~A.~Athanasiadis}
\address{Department of Mathematics\\
National and Kapodistrian University of Athens\\
Panepistimioupolis\\ 15784 Athens, Greece}
\email{caath@math.uoa.gr}

\date{January 24, 2026}

\begin{abstract}
An $n$-dimensional lattice polytope $\qQ_\sigma$ 
can be associated to any composition $\sigma$ of 
a positive integer $n$, as a special case of 
constructions due to Pitman--Stanley and Chapoton. 
The entries of the $h$-vector of $\sigma$, 
introduced by Chapoton, enumerate the lattice 
points in $\qQ_\sigma$ by the number of their 
nonzero coordinates. Chapoton conjectured that 
this vector is equal to the $h$-vector of a 
flag simplicial polytope. This paper proves this 
conjecture. Moreover, it shows that the 
gamma-vector associated to the $h$-vector of 
$\sigma$ is nonnegative by means of an explicit 
combinatorial interpretation and confirms certain 
other conjectures of Chapoton on the lattice point 
enumeration of composition polytopes. A 
combinatorial interpretation of their 
$h^\ast$-polynomials is deduced. 
\end{abstract}

\maketitle

\section{Introduction}
\label{sec:intro}
 
A composition of a positive integer $n$ is 
defined as a sequence $\sigma = (r_1, r_2,\dots,r_k)$ 
of positive integers whose sum is equal to $n$. 
Given such a composition $\sigma$, we consider the 
polytope $\qQ_\sigma$ in $\RR^n$ defined by the 
inequalities $x_i \ge 0$ for $1 \le i \le n$ and 
\begin{equation} \label{eq:Q-def}
x_1 + x_2 + \cdots + x_{s_i} \le s_i
\end{equation}
for $1 \le i \le k$, where $s_i = r_1 + r_2 + 
\cdots + r_i$ and $x_1, x_2,\dots,x_n$ are the 
standard coordinates of $\RR^n$. We refer to the 
polytopes $\qQ_\sigma$ in this paper as 
\emph{composition polytopes}.

Composition polytopes are special cases of the 
polytopes associated to arbors (defined as rooted 
trees on the vertex set of blocks of a partition 
of a finite set), introduced and studied by 
Chapoton~\cite{Cha25+}. They also arise as 
Pitman--Stanley polytopes~\cite{SP02} when some of 
the parameters defining these polytopes are 
allowed to be zero. The composition polytope 
$\qQ_\sigma$ reduces to the $n$th dilate of the 
$n$-dimensional simplex in $\RR^n$ with vertices 
the origin and the unit coordinate vectors, when 
$\sigma$ has one part, and to the Pitman--Stanley 
polytope with all parameters equal to 1, when all 
parts of $\sigma$ are equal to 1.  

The face lattice of composition polytopes is well
understood. As observed in 
\cite[Theorem~20]{SP02}, $\qQ_\sigma$ is 
combinatorially isomorphic to a product of 
simplices, the dimensions of which are the parts 
of $\sigma$. The lattice point enumeration of 
composition polytopes, which is more interesting, 
has been studied for the more general classes of 
arbor polytopes~\cite{Cha25+} and Pitman--Stanley 
polytopes~\cite{SP02}. Chapoton defined the 
\emph{$h$-vector} of a composition $\sigma$ of $n$ 
(more generally, of an arbor) as the sequence 
$h(\sigma) = (h_i(\sigma))_{0 \le i \le n}$, where 
$h_i(\sigma)$ is the number of points in $\qQ_\sigma 
\cap \NN^n$ having exactly $i$ nonzero coordinates
(this vector should not be confused with the 
$h$-vector of the simple polytope $\qQ_\sigma$). 
We note that $h_0(\sigma) = h_n(\sigma) = 1$. 
The following statement is a special case of an 
intriguing conjecture of Chapoton 
\cite[Conjecture~0.1]{Cha25+}, stated for all 
arbor polytopes; we refer to the discussion in 
\cite[Section~0.1]{Cha25+} for the motivation 
behind it.
\begin{conjecture} [\cite{Cha25+}] \label{conj:poly}
The vector $h(\sigma)$ is equal to the $h$-vector 
of an $n$-dimensional simplicial polytope for every 
composition $\sigma$ of $n$.
\end{conjecture}

In particular, $h(\sigma)$ should be palindromic
and unimodal, meaning that $h_i(\sigma) = h_{n-i}
(\sigma)$ for all $0 \le i \le n$ and $h_0(\sigma) 
\le h_1(\sigma) \le \cdots \le 
h_{\lfloor n/2\rfloor}(\sigma)$. Moreover, it is 
suggested in \cite[Section~0.1]{Cha25+} that 
$h(\sigma)$ should in fact be equal to the 
$h$-vector of an $n$-dimensional flag simplicial 
polytope (a simplicial polytope $\qQ$ is said to 
be \emph{flag} if every clique in the one-skeleton 
of $\qQ$ is the vertex set of a face of $\qQ$). 
Thus, in view of Gal's conjecture 
\cite[Conjecture~2.1.7]{Ga05} (see also 
\cite[Section~3]{Ath18}) on gamma-positivity, it 
is reasonable to expect that the generating 
polynomial $h(\sigma,t) := \sum_{i=0}^n 
h_i(\sigma) t^i$ of $h(\sigma)$ is 
$\gamma$-positive, meaning that 
\begin{equation} \label{eq:tau-gamma}
h(\sigma,t) = \sum_{i=0}^{\lfloor n/2 \rfloor} 
\gamma_i(\sigma) t^i (1+t)^{n-2i}
\end{equation}
for some nonnegative integers $\gamma_i(\sigma)$. 
The first main result of this paper states that 
all these claims are true for composition polytopes.
\begin{theorem} \label{thm:gamma-tau} 
The vector $h(\sigma)$ is equal to the $h$-vector 
of an $n$-dimensional flag simplicial polytope 
for every composition $\sigma$ of $n$. Moreover, 
$h(\sigma,t)$ is $\gamma$-positive and, in 
particular, palindromic and unimodal.
\end{theorem}

Our second main result confirms another curious
conjecture of Chapoton~\cite{Cha25+}, concerning 
the Ehrhart polynomial of composition polytopes.
We recall (see, for instance, 
\cite[Section~4.6]{StaEC1}) that the \emph{Ehrhart 
polynomial} $\Ehr(\qQ, t)$ of a lattice polytope 
$\qQ \subset \RR_{\ge 0}^n$ has the defining 
property that $\Ehr(\qQ, m) = \# (m\qQ \cap 
\NN^n)$ for every $m \in \NN$. The \emph{zeta 
polynomial} (see~\cite[Section~3.12]{StaEC1}) 
$\zZ(P, t)$ of a finite poset $(P, \preceq)$ 
has the defining property that $\zZ(P, m+1)$ is 
equal to the number of $m$-element multichains 
in $P$, meaning sequences $(p_1, p_2,\dots,p_m)$ 
of elements of $P$ such that $p_1 \preceq p_2 
\preceq \cdots \preceq p_m$, for every $m \in 
\NN$. Given a composition $\sigma$ of $n$, 
Chapoton~\cite[Section~1]{Cha25+} has considered
the set $P_\sigma$ of lattice points in $\qQ_\sigma$ 
as a subposet of the product poset $\NN^n$.

\begin{theorem} \label{thm:EZ} 
{\rm (\cite[Conjecture~3.1]{Cha25+})} We have 
$\Ehr(\qQ_\sigma, t) = \zZ(P_{\rev(\sigma)}, t+1)$ 
for every composition $\sigma$, where $\rev(\sigma)$ 
is the composition obtained from $\sigma$ by 
reversing the order in which the parts of $\sigma$ 
appear.
\end{theorem}

Here is another intriguing conjecture 
from \cite{Cha25+}, stated there for all arbor 
polytopes (see \cite[Conjecture~2.2]{Cha25+}), 
on which no progress is made in this paper.
\begin{conjecture} [\cite{Cha25+}]
For every composition $\sigma$, all roots of the 
Ehrhart polynomial of $\qQ_\sigma$ are real and 
lie in the interval $[-1,0]$.
\end{conjecture}

Theorem~\ref{thm:gamma-tau} is proven in 
Section~\ref{sec:polytope}. This section shows 
that $h(\sigma,t)$ is equal to the Eulerian 
polynomial of a graded poset of width two (see 
Theorem~\ref{thm:poset-tau}) and deduces 
Theorem~\ref{thm:gamma-tau} from results of 
Br\"and\'en~\cite{Bra04} and 
Reiner--Welker~\cite{RW05} on poset Eulerian 
polynomials. Section~\ref{sec:gamma} provides a 
direct proof of the gamma-positivity of 
$h(\sigma,t)$ via an explicit combinatorial 
interpretation of the integers $\gamma_i(\sigma)$ 
which appear in Equation~(\ref{eq:tau-gamma}) 
(see Proposition~\ref{prop:gamma}); it also 
shows that $h(\sigma) = h(\rev(\sigma))$ for 
every composition $\sigma$ (see 
Proposition~\ref{prop:rev-h}). Two proofs of 
Theorem~\ref{thm:EZ}, one essentially bijective 
and another based on the results of 
\cite[Section~3]{SP02}, are given in 
Section~\ref{sec:EZ} and a combinatorial 
interpretation of the $h^\ast$-polynomial of 
$\qQ_\sigma$ is deduced via the theory of pure 
shellability of posets (see 
Proposition~\ref{prop:h-star}). While this paper 
is restricted to composition polytopes, we do 
hope that some of the arguments presented here
can be extended to more general classes of 
lattice polytopes.

Throughout this paper,   
we set $\NN = \{0, 1, 2,\dots\}$ and $[n] = 
\{1, 2,\dots,n\}$ for $n \in \NN$. We denote by 
$\# S$ the cardinality of a finite set $S$ and 
write $\sigma \models n$ when $\sigma$ is a 
composition of $n$. Our notation and terminology 
about posets follows \cite[Chapter~3]{StaEC1}. 
Background on Ehrhart theory and face enumeration 
of simplicial complexes and polytopes can be found 
in~\cite{BR15,HiAC,StaCCA} 
\cite[Chapter~4]{StaEC1}. 

\section{Gamma-positivity}
\label{sec:gamma}

This section provides a direct combinatorial 
proof of the gamma-positivity of $h(\sigma,t)$ 
and shows that $h(\sigma) = h(\rev(\sigma))$ for 
every composition $\sigma \models n$. Moreover, it 
introduces a lattice path interpretation of points 
in $P_\sigma = \qQ_\sigma \cap \NN^n$ which will 
be useful in the following section as well.

Throughout this section we let $\sigma = (r_1, 
r_2,\dots,r_k)$ be a composition of $n$ and set 
$s_i = r_1 + r_2 + \cdots + r_i$ for $i \in [k]$. 
We recall that $P_\sigma$ consists of all points 
$(a_1, a_2,\dots,a_n) \in \NN^n$ such that
\[ a_1 + a_2 + \cdots + a_{s_i} \le s_i \]

\medskip
\noindent
for every $i \in [k]$. We record such a point by 
the $n$-tuple $(p_1, p_2,\dots,p_n) \in \NN^n$, 
where 
\[ p_i = \begin{cases}
         a_1 + a_2 + \cdots + a_i, & 
				 \text{if $a_i \ge 1$} \\
				 0, & \text{if $a_i = 0$} \\
         \end{cases} \]
for every $i \in [n]$. As a running example, if 
$\sigma = (3, 4, 2)$, then $(0, 0, 2, 0, 4, 0, 0, 
3, 0) \in P_\sigma$ is recorded by the tuple $(0, 0, 
2, 0, 6, 0, 0, 9, 0)$. The map just described shows 
that $P_\sigma$ is in one-to-one correspondence 
with the set, say $\hH_\sigma$, of $n$-tuples $(p_1, 
p_2,\dots,p_n) \in \NN^n$ consisting of positive 
integers which appear in strictly increasing order 
from left to right and zeros, such that $p_1, 
p_2,\dots,p_{s_i} \le s_i$ for every $i \in [k]$. 
Moreover, $h_j(\sigma)$ is equal to the number of 
elements of $\hH_\sigma$ which have exactly $j$ 
nonzero entries.  

Given $x \in \hH_\sigma$, we denote by $S(x)$ the 
set of positions of the nonzero entries in $x$
and by $T(x)$ the set of these entries. Thus, in 
our running example, we have $S(x) = \{3, 5, 8\}$
and $T(x) = \{2, 6, 9\}$. Clearly, $S(x), T(x)
\subseteq [n]$ and $\# S(x) = \# T(x)$ for every 
$x \in \hH_\sigma$. Moreover, every $x \in 
\hH_\sigma$ is determined by $S(x)$ and $T(x)$
and 

\begin{equation} \label{eq:tau-H}
h(\sigma,t) = \sum_{x \in \hH_\sigma} t^{\# S(x)}.
\end{equation}

\smallskip
\begin{lemma} \label{lem:gamma} 
Let $S, T$ be subsets of $[n]$ of equal 
cardinalities and let $R \subseteq [n]$ be disjoint
from $S \cup T$. Then, the pair $(S,T)$ has the 
property that $S = S(x)$ and $T = T(x)$ for some 
$x \in \hH_\sigma$ if and only if so does $(S \cup 
R, T \cup R)$.
\end{lemma}

\begin{proof}
We may assume that $R$ has a single element $r \in 
[n]$. Let $x = (p_1, p_2,\dots,p_n) \in \NN^n$ be 
the $n$-tuple which has the elements of $T$ as its 
nonzero entries and these appear in strictly 
increasing order from left to right in the 
positions of the elements of $S$. Similarly, let 
$y = (q_1, q_2,\dots,q_n) \in \NN^n$ be the 
$n$-tuple which has the elements of $T \cup \{r\}$ 
as its nonzero entries and these appear in strictly 
increasing order from left to right in the positions 
of the elements of $S \cup \{r\}$. We have $x \in 
\hH_\sigma$ if and only if 
\begin{equation} \label{eq:p-cond}
p_1, p_2,\dots,p_{s_i} \le s_i \ \ \text{for every} 
\ i \in [k]
\end{equation}
and 
$y \in \hH_\sigma$ if and only if 
\begin{equation} \label{eq:q-cond}
q_1, q_2,\dots,q_{s_i} \le s_i \ \ \text{for every} 
\ i \in [k].
\end{equation}
Thus, we have to show the equivalence of 
(\ref{eq:p-cond}) and (\ref{eq:q-cond}). This is 
obvious if $q_r = r$, since then $x$ and $y$ differ
only in the $r$th coordinate. Suppose that $r = q_s$
for some $s>r$. Then, $y$ is obtained from $x$ by 
placing $r$ in position $s$ and shifting to the 
left some nonzero entries, all less than $r$, 
within positions $r, r+1,\dots,s$. Finally, if $r = 
q_s$ for some $s<r$, then $y$ is obtained from $x$ 
by placing $r$ in position $s$ and shifting to the 
right some nonzero entries, all larger than $r$, 
within positions $s, s+1,\dots,r$. The equivalence
(\ref{eq:p-cond}) $\Leftrightarrow$ 
(\ref{eq:q-cond}) holds in either case.
\end{proof}

Theorem~\ref{thm:gamma-tau} is a consequence of 
the following result.

\begin{proposition} \label{prop:gamma} 
Equation~(\ref{eq:tau-gamma}) holds if $\gamma_i
(\sigma)$ is defined as the number of elements $x 
\in \hH_\sigma$ such that $S(x) \cap T(x) = 
\varnothing$ and $\# S(x) = i$, for $0 \le i \le 
\lfloor n/2 \rfloor$.
\end{proposition}

\begin{proof}
For $x, y \in \hH_\sigma$ let us write 
$x \sim_\sigma y$ if 
\begin{align*}
S(x) \sm (S(x) \cap T(x)) & = S(y) \sm (S(y) \cap 
T(y)) \\
T(x) \sm (S(x) \cap T(x)) & = T(y) \sm (S(y) \cap 
T(y)).
\end{align*}
Clearly, $\sim_\sigma$ is an 
equivalence relation on $\hH_\sigma$. In view of 
Lemma~\ref{lem:gamma}, every equivalence class has
a unique representative $x \in \hH_\sigma$ such 
that $S(x) \cap T(x) = \varnothing$ and the map 
which assigns the pair $(S(y), T(y))$ to $y \in 
\hH_\sigma$ induces a bijection from the 
equivalence class, say $C_x$, of $x$ to the set 
of pairs $(S(x) \cup R, T(x) \cup R)$, where $R$
ranges over all subsets of $[n] \sm (S(x) \cup 
T(x))$. As a result, 

\[ \sum_{y \in C_x} t^{\# S(y)} = 
   \sum_{R \subseteq [n] \sm (S(x) \cup T(x))} 
   t^{\# (S(x) \cup R)} = t^{\# S(x)} 
	 (1+t)^{n - 2 \, \# S(x)}. \]

\medskip
\noindent
We sum over all equivalence classes of 
$\sim_\sigma$ and the proof follows 
from~(\ref{eq:tau-H}).
\end{proof}

\begin{example} \label{ex:NarB} \rm 
Let $\sigma = (n)$ be the composition of $n$ with
only one part. Then, $P_\sigma$ consists of all 
points $(a_1, a_2,\dots,a_n) \in \NN^n$ such that
$a_1 + a_2 + \cdots + a_n \le n$, $h_i(\sigma)$
is equal to the total number of pairs $(S,T)$ of 
$i$-element subsets of $[n]$ and $\gamma_i(\sigma)$
is equal to the number of such pairs $(S,T)$ with 
$S \cap T = \varnothing$. Hence, $h_i(\sigma) = 
{n \choose i}^2$ for $0 \le i \le n$ and 
$\gamma_i(\sigma) = {n \choose i}{n-i \choose i}$
for $0 \le i \le \lfloor n/2 \rfloor$ and 
Proposition~\ref{prop:gamma} specializes to the
well known identity 
\[ \sum_{i=0}^n {n \choose i}^2 t^i = 
   \sum_{i=0}^{\lfloor n/2 \rfloor} 
	 {n \choose i}{n-i \choose i} t^i (1+t)^{n-2i} 
	 \]
(see \cite[Theorem~2.32]{Ath18} for a 
generalization in the setting of finite Coxeter 
groups).
\qed
\end{example}

\noindent
\textbf{Lattice path interpretation}. For the 
remainder of this section we represent elements 
of $\NN^n$ with sum of coordinates at most $n$ 
by lattice paths in the $xy$-plane as follows. 
A point $(a_1, a_2,\dots,a_n) \in \NN^n$ with 
$a_1 + a_2 + \cdots + a_n \le n$ is represented 
by the lattice path with initial vertex $(0,0)$, 
terminal vertex $(n, n)$ and east (E) or north 
(N) unit steps, so that the height of the part 
of the lattice path within the strip $j-1 \le x 
\le j$ is $a_1 + a_2 + \cdots + a_j$ for every 
$j \in [n]$. We observe that the number of 
nonzero entries of such an element of $\NN^n$ 
is equal to the number of north-east corners 
of the associated lattice path. 

We denote by $\gG_n$ the set of lattice 
paths which represent the elements of $\NN^n$
with sum of coordinates at most $n$ and let 
$\Gamma_\sigma \in \gG_n$ be the lattice path 
associated to the $n$-tuple with entry $s_i$ 
in position $s_{i-1} + 1$ for $1 \le i \le k$ 
(where $s_0 := 0$) and all other entries equal 
to zero. Thus, $\Gamma_\sigma$ has $r_1$ north 
steps, followed by $r_1$ east steps, followed 
by $r_2$ north steps, followed by $r_2$ east 
steps, and so on. Then, the elements of 
$P_\sigma$ are represented by the lattice 
paths in $\gG_n$ which lie weakly below 
$\Gamma_\sigma$. 
\begin{example} \label{ex:NarA} \rm 
Let $\sigma = (1, 1,\dots,1) \models n$, so that 
$P_\sigma$ consists of all points $(a_1, 
a_2,\dots,a_n) \in \NN^n$ such that $a_1 + a_2 + 
\cdots + a_j \le j$ for every $j \in [n]$. The 
lattice path $\Gamma_\sigma$ has successive 
steps $\textrm{NENE} \cdots \textrm{NE}$ and the
total number of lattice paths in $\gG_n$ which 
lie weakly below it (and hence the cardinality 
of $P_\sigma$) is well known to equal the 
Catalan number $\frac{1} {n+2}{2n+2 \choose n+1}$. 
Moreover, counting these lattice paths by the 
number of north-east corners shows that $h_i
(\sigma)$ is equal to the Narayana number
\[ h_i(\sigma) = \frac{1}{i+1} {n \choose i} 
   {n+1 \choose i} \]

\medskip
\noindent
for $0 \le i \le n$ (see, for instance, 
\cite[Problem~A.46]{StaCat}). 

By Proposition~\ref{prop:gamma}, $\gamma_i
(\sigma)$ is equal to the number of pairs $(S,T)$ 
of disjoint $i$-element subsets of $[n]$ such that 
the $r$th largest element of $T$ is smaller that
the $r$th largest element of $S$ for every $r$. 
Since these pairs are in one-to-one 
correspondence with Young tableaux of shape 
$(i,i)$ with entries in $[n]$, we have 
\begin{align*} 
\gamma_i(\sigma) & = {n \choose 2i} \cdot 
\frac{1}{i+1} {2i \choose i} = \frac{1}{i+1} 
\cdot \frac{n!}{(i!)^2 \, (n-2i)!} \\ & 
= \frac{1}{i+1}{n \choose i}{n-i \choose i} 
\end{align*}

\medskip
\noindent
for $0 \le i \le \lfloor n/2 \rfloor$ and 
Proposition~\ref{prop:gamma} specializes to the
well known identity 
\[ \sum_{i=0}^n \frac{1}{i+1} {n \choose i} 
   {n+1 \choose i} t^i = 
	 \sum_{i=0}^{\lfloor n/2 \rfloor} \frac{1}{i+1}
	 {n \choose i}{n-i \choose i} t^i (1+t)^{n-2i} 
	 \]
(see, for instance, again 
\cite[Theorem~2.32]{Ath18}).
\qed
\end{example}

We conclude this section with the following 
statement.

\begin{proposition} \label{prop:rev-h} 
{\rm (\cite[Conjecture~3.2]{Cha25+})} We have 
$h(\sigma) = h(\rev(\sigma))$ for every $\sigma 
\models n$.
\end{proposition}

\begin{proof}
We observe that orthogonal reflection in the 
line $x+y = n$ maps $\Gamma_\sigma$ to 
$\Gamma_{\rev(\sigma)}$ and hence induces a 
bijection from the set of lattice paths which 
represent the elements of $P_\sigma$ to the set 
of those for $P_{\rev(\sigma)}$. We also observe
that the number of north-east corners of a 
lattice path is preserved under this bijection 
and the proof follows.
\end{proof}

\section{Polytopality}
\label{sec:polytope}

This section proves Theorem~\ref{thm:gamma-tau}.
We first need to recall some background and 
terminology about the face enumeration of 
simplicial polytopes and the combinatorics of 
posets; more information can be found in 
\cite{RW05, StaCCA} \cite[Section~3.15]{StaEC1}.

The \emph{$h$-polynomial} of a $d$-dimensional 
simplicial polytope $\qQ$ is a convenient way
to record the face numbers of $\qQ$; it is 
defined by the formula 
\begin{equation} \label{eq:h-poly}
h_\qQ(t) = \sum_{i=0}^d f_{i-1}(\qQ) t^i 
           (1-t)^{d-i},
\end{equation}
where $f_{i-1}(\qQ)$ is the number of 
$(i-1)$-dimensional faces of $\qQ$. If $h_\qQ
(t) = \sum_{i=0}^d h_i(\qQ) t^i$, then 
$(h_0(\qQ), h_1(\qQ),\dots,h_d(\qQ))$ is called
the \emph{$h$-vector} of $\qQ$. The polynomial 
$h_\qQ(t)$ has nonnegative coefficients and, as 
a consequence of the $g$-theorem 
\cite[Section~III.1]{StaCCA}, it is palindromic 
and unimodal for every simplicial polytope $\qQ$.
Moreover, it is conjectured to be 
$\gamma$-positive \cite[Conjecture~2.1.7]{Ga05} 
whenever $\qQ$ is flag.

Let $P$ be a poset with $p$ elements, with  
partial order denoted by $\preceq$, and 
$\omega: P \to [p]$ be a natural labeling,
meaning a bijection with the property that $a 
\prec b \Rightarrow \omega(a) < \omega(b)$ for 
all $a, b \in P$. A \emph{linear extension} of
$P$ is defined as any permutation $(a_1, 
a_2,\dots,a_p)$ of $P$ such that $a_i \prec a_j
\Rightarrow i < j$ for all $i, j \in [p]$. The
\emph{Jordan-H\"older set} $\lL(P,\omega)$ 
consists of all permutations $(\omega(a_1), 
\omega(a_2),\dots,\omega(a_p))$ of $[p]$, where
$(a_1, a_2,\dots,a_p)$ is a linear extension of 
$P$. The \emph{Eulerian polynomial} of $P$ is 
defined as 
\begin{equation} \label{eq:P-Eulerian}
W_P(t) = \sum_{w \in \lL(P,\omega)} t^{\des(w)},
\end{equation}
where $\des(w)$ stands for the number of 
indices $i \in [p-1]$ (called \emph{descents}) 
such that $w_i > w_{i+1}$,
if $w = (w_1, w_2,\dots,w_p)$. This polynomial 
is independent of the natural labeling $\omega$;
it was shown to be $\gamma$-positive and equal 
to the $h$-polynomial of a $(p-d)$-dimensional 
simplicial polytope by Br\"and\'en~\cite{Bra04} 
and Reiner--Welker~\cite{RW05}, respectively, 
for every graded poset $P$, where $d$ is the 
cardinality of any maximal chain of $P$.

The \emph{width} of a poset $P$ is 
defined as the largest size of an antichain 
(subset of pairwise incomparable elements) in 
$P$. The following statement is a stronger 
version of Theorem~\ref{thm:gamma-tau}.
\begin{theorem} \label{thm:poset-tau} 
The polynomial $h(\sigma,t)$ is equal to the 
Eulerian polynomial of a graded poset of width 
two for every composition $\sigma \models n$. 

In particular, $h(\sigma,t)$ is $\gamma$-positive,
hence palindromic and unimodal, and equal to the 
$h$-polynomial of an $n$-dimensional flag 
simplicial polytope.
\end{theorem}

\begin{proof}
Let $\sigma = (r_1, r_2,\dots,r_k) \models n$ 
and $s_i = r_1 + r_2 + \cdots + r_i$ for $i \in 
[k]$ be as in Section~\ref{sec:gamma}. We 
consider a poset $\Pi_\sigma$ which consists of 
two disjoint chains $\alpha_1 \prec \alpha_2 
\prec \cdots \prec \alpha_n$ and $\beta_1 \prec 
\beta_2 \prec \cdots \prec \beta_n$, together 
with the additional cover relations 
$\alpha_{s_i} \prec \beta_{s_i+1}$ for every 
$i \in [k-1]$. This poset is graded; all maximal 
chains have cardinality $n$. We also consider the 
natural labeling $\omega: \Pi_\sigma \to [2n]$ 
such that $\omega(\alpha_i) = i$ and 
$\omega(\beta_i) = n+i$ for $i \in [n]$.

We claim that $h(\sigma,t) = W_{\Pi_\sigma}(t)$.
Indeed, we recall from Section~\ref{sec:gamma}
that the coefficient $h_j(\sigma)$ of $t^j$ in 
$h(\sigma,t)$ is equal to the number of lattice
paths in $\gG_n$ which lie weakly below 
$\Gamma_\sigma$ and have exactly $j$ north-east
corners. On the other hand, the coefficient of 
$t^j$ in $W_{\Pi_\sigma}(t)$ is equal to the 
number of permutations $w \in \lL(\Pi_\sigma,
\omega)$ such that $\des(w) = j$. To any linear
extension of $\Pi_\sigma$ one can associate the
lattice path in $\gG_n$ obtained by replacing
each element $\alpha_i$ by an east step and 
each element $\beta_i$ by a north step. These 
lattice paths lie weakly below $\Gamma_\sigma$
and the resulting map induces a one-to-one 
correspondence between permutations $w \in \lL
(\Pi_\sigma,\omega)$ and lattice paths in $\gG_n$ 
which lie weakly below $\Gamma_\sigma$. Moreover, 
the number of descents of $w \in \lL(\Pi_\sigma,
\omega)$ is equal to the number of north-east 
corners of the corresponding lattice path. These 
facts show that $[t^j] W_{\Pi_\sigma}
(t) = h_j(\sigma)$ for every $j$ and the proof 
of the claim follows. The first statement in 
the theorem follows from the claim. The second 
statement follows from the first and the 
aforementioned results of \cite{Bra04, RW05}.  
Since $\Pi_\sigma$ has width two, the simplicial
polytope with $h$-polynomial $h(\sigma,t)$ can be
chosen to be flag by \cite[Theorem~3.23]{RW05}.
\end{proof}

\begin{remark} \label{rem:Cha} \rm
After this paper first appeared on the arXiv,
Frederic Chapoton~\cite{Cha25} suggested to the 
author that $h(\sigma,t)$ may be equal to the 
$h^\ast$-polynomial of the order polytope of the
poset $\Pi_\sigma$. This suggestion led to 
Theorem~\ref{thm:poset-tau}. 

The lattice path interpretation of $h(\sigma,t)$
shows that it belongs to the class of nonnesting
rook polynomials, introduced and studied recently
by Alexandersson and Jal~\cite{AJ24+}. Given that, 
Theorem~\ref{thm:poset-tau} can also be deduced 
from the results of \cite[Section~5]{AJ24+}. 
Corollary~4.2 in~\cite{AJ24+} implies that 
$h(\sigma,t)$ has another important property, 
namely that it is ultra log-concave for every 
composition $\sigma$. To the best of our 
knowledge, the real-rootedness of $h(\sigma,t)$ 
is open.
\end{remark}

\section{Ehrhart and Zeta polynomials}
\label{sec:EZ}

This section offers two proofs of 
Theorem~\ref{thm:gamma-tau} and deduces a 
combinatorial interpretation of the 
$h^\ast$-polynomial of the composition polytope 
$\qQ_\sigma$. As in the previous sections, $\sigma 
= (r_1, r_2,\dots,r_k)$ is a composition of $n$ 
and $s_i = r_1 + r_2 + \cdots + r_i$ for $0 \le 
i \le k$, where $s_0 := 0$.

We recall that $\NN^n$ is partially ordered as 
a direct product, meaning that for $x, y \in 
\NN^n$ with $x = (a_1, a_2,\dots,a_n)$ and $y = 
(b_1, b_2,\dots,b_n)$ we have $x \preceq y$ if 
and only if $a_i \le b_i$ for every $i \in [n]$,
where $\preceq$ stands for the partial order of 
$\NN^n$.
\begin{lemma} \label{lem:multichains} 
There is a one-to-one correspondence from the 
set of $m$-multichains in $\NN^n$ with 
a given maximal element $(a_1, a_2,\dots,a_n)$ 
to the set of solutions 
$(b_{i,j})_{1 \le i \le n, \, 1 \le j \le m} \in 
\NN^{mn}$ of the system of linear equations 
\[ x_{i,1} + x_{i,2} + \cdots + x_{i,m} = a_i \]
for $i \in [n]$.
\end{lemma}

\begin{proof}
By definition of the product order $\NN^n$, it 
suffices to verify the statement for $n=1$, which
is a trivial task.
\end{proof} 

\begin{proof}[First proof of Theorem~\ref{thm:EZ}]
For any $m \in \NN$, $\Ehr(\qQ_\sigma,m)$ is
equal to the number of points $(a_1, a_2,\dots,a_n)
\in \NN^n$ such that 
\begin{equation} \label{eq:ineq-ai}
a_1 + a_2 + \cdots + a_{s_i} \le s_i m
\end{equation}
for $i \in [k]$. Thus, for $i=k$ we get $a_1 + 
a_2 + \cdots + a_n \le mn$. Let's consider the 
partial sums $a_1 + a_2 + \cdots + a_i$ and, for 
$0 \le j \le mn$, let's denote by $b_j$ the number 
of indices $i \in [n]$ for which $a_1 + a_2 + 
\cdots + a_i = j$ (the reader may want to 
interpret the $b_j$ when $(a_1, a_2,\dots,a_n)$ 
is represented by a lattice path, as in 
Section~\ref{sec:gamma}). Then, 
\begin{equation} \label{eq:sum-bj}
b_0 + b_1 + \cdots + b_{mn} = n
\end{equation}
and 
\begin{equation} \label{eq:ineq-bj-1}
b_0 + b_1 + \cdots + b_{s_i m} \ge s_i
\end{equation}
for $i \in [k]$. Moreover, the map which sends
$(a_1, a_2,\dots,a_n)$ to $(b_0, b_1,\dots,b_{mn})$
is a one-to-one correspondence from the set of 
solutions $(a_1, a_2,\dots,a_n) \in \NN^n$ of the 
system of linear inequalities (\ref{eq:ineq-ai})
to the set of solutions $(b_0, b_1,\dots,b_{mn})
\in \NN^{mn+1}$ of the system of linear 
equalities and inequalities (\ref{eq:sum-bj}) and 
(\ref{eq:ineq-bj-1}). Eliminating $b_0$ we find 
that the latter set is in one-to-one correspondence 
with the set of solutions $(b_1, b_2,\dots,b_{mn}) 
\in \NN^{mn}$ of the system of linear inequalities 
\begin{equation} \label{eq:ineq-bj-2}
b_{s_i m + 1} + b_{s_i m + 2} + \cdots + b_{mn} 
       \le n - s_i
\end{equation}
for $i \in [k-1]$ and 
\begin{equation} \label{eq:ineq-bj}
b_1 + b_2 + \cdots + b_{mn} \le n.
\end{equation}
The number of these solutions is equal to the 
number $\zZ(P_{\rev(\sigma)},m+1)$ of 
$m$-multichains in $P_{\rev(\sigma)}$. This 
follows from the fact that, since $n - s_i = 
r_{i+1} + r_{i+2} + \cdots + r_k$, the set 
$P_{\rev(\sigma)}$ consists of all tuples 
$(x_1, x_2,\dots,x_n) \in \NN^n$ such that
\[ x_{s_i + 1} + x_{s_i + 2} + \cdots + x_n 
       \le n - s_i \]
for $0 \le i < k$ and from
Lemma~\ref{lem:multichains}.
\end{proof}

Part of our second proof of Theorem~\ref{thm:EZ} 
follows the proof of \cite[Equation~(33)]{SP02},
which is a formula for the number of lattice 
points of a Pitman--Stanley polytope. To avoid 
making too long of a digression, we assume 
familiarity with Sections~3 and~4 of \cite{SP02} 
and only sketch this part of the proof.
\begin{proof}
[Second proof of Theorem~\ref{thm:EZ}]
Let us consider the poset $\Omega_\sigma$ which
consists of two disjoint chains $\alpha_1 \prec
\alpha_2 \prec \cdots \prec \alpha_n$ and 
$\beta_1 \prec \beta_2 \prec \cdots \prec 
\beta_n$, together with the additional cover 
relations $\alpha_{s_i} \prec \beta_{s_i}$ for 
every $i \in [k]$. We let $C = \{\beta_1, 
\beta_2,\dots,\beta_n\}$ and $\bu = (u_1, 
u_2,\dots,u_n)$, where $u_j = jm$ for every
$j \in [n]$ and where $m \in \NN$ is fixed. In 
this situation, Pitman and 
Stanley~\cite[Section~3]{SP02} have defined the 
order cone section $\cC_C(\Omega_\sigma,\bu)$ as 
the polytope of all order preserving maps $f: 
\Omega_\sigma \to \RR_{\ge 0}$ such that 
$f(\beta_j) = u_j$ for every $j \in [n]$. 
Setting $f_j = f(\alpha_j)$ for every $j \in 
[n]$ and $f_0 = 0$ we have 
\[ 0 = f_0 \le f_1 \le \cdots \le f_n \]
and 
\[ f_{s_i} \le s_i m \] 
for every $i \in [k]$. The transformation $x_j 
= f_j - f_{j-1}$ for $j \in [n]$ then shows 
that $\cC_C(\Omega_\sigma,\bu)$ is integrally 
equivalent to $m \qQ_\sigma$. Theorem~5 in 
\cite{SP02} expresses the number of lattice 
points of $\cC_C(\Omega_\sigma,\bu)$ as a sum
over all linear extensions of $\Omega_\sigma$.
Such a linear extension is uniquely determined
by the positions, say $1 \le j_1 < j_2 < \cdots 
< j_n \le 2n$, of $\beta_1, 
\beta_2,\dots,\beta_n$, which must satisfy 
$j_{s_i} \ge 2 s_i$ for every $i \in [k]$. 
Thus, setting $\kappa_i = j_i - j_{i-1} - 1$ 
for $i \in [n]$ (where $j_0 = 0$), the argument 
of \cite[p.~617]{SP02} shows that in our 
situation, \cite[Theorem~5]{SP02} yields the
formula 
\begin{equation} \label{eq:SP} 
\Ehr(\qQ_\sigma, m) = 
\sum_{(\kappa_1, \kappa_2,\dots,\kappa_n) \in 
K_\sigma} {m + \kappa_1 \choose \kappa_1} 
{m + \kappa_2 - 1 \choose \kappa_2} \cdots 
{m + \kappa_n - 1 \choose \kappa_n}, 
\end{equation}
where $K_\sigma$ consists of all tuples $(\kappa_1, 
\kappa_2,\dots,\kappa_n) \in \NN^n$ such that
$\kappa_1 + \kappa_2 + \cdots + \kappa_n = n$
and $\kappa_1 + \kappa_2 + \cdots + \kappa_{s_i} 
\ge s_i$ for every $i \in [k]$.

Thus, it suffices to prove that the right-hand 
side of Equation~(\ref{eq:SP}) is also equal to 
the number $\zZ(P_{\rev(\sigma)}, m+1)$ of 
$m$-element multichains in $P_{\rev(\sigma)}$. 
As already observed, 
$P_{\rev(\sigma)}$ consists of all tuples 
$(\lambda_1, \lambda_2,\dots,\lambda_n) \in 
\NN^n$ such that $\lambda_1 + \lambda_2 + \cdots
+ \lambda_n \le n$ and
\[ \lambda_{s_i + 1} + \lambda_{s_i + 2} + 
   \cdots + \lambda_n \le n - s_i \]
for $i \in [k-1]$. Setting $\lambda_0 = n - 
\lambda_1 - \lambda_2 - \cdots - \lambda_n \ge 0$, 
we have $\lambda_0 + \lambda_1 + \cdots + 
\lambda_n = n$ and 
\[ \lambda_0 + \lambda_1 + \cdots + \lambda_{s_i} 
   \ge s_i \]
for every $i \in [k]$ and $\zZ(P_{\rev(\sigma)}, 
m+1)$ is equal to the sum of the number of 
$m$-multichains in $\NN^n$ with maximum element 
$(\lambda_1, \lambda_2,\dots,\lambda_n)$ over all 
such tuples $(\lambda_0, 
\lambda_1,\dots,\lambda_n) \in \NN^{n+1}$. We may
set $\kappa_1 = \lambda_0 + \lambda_1$ and 
$\kappa_i = \lambda_i$ for all $1 < i \le n$ and
sum instead the number of $m$-multichains in 
$\NN^n$ with maximum element $(\lambda_1, 
\kappa_2,\dots,\kappa_n)$, where $0 \le \lambda_1 
\le \kappa_1$, over all tuples $(\kappa_1, 
\kappa_2,\dots,\kappa_n) \in K_\sigma$. By 
Lemma~\ref{lem:multichains} and basic enumerative
combinatorics, the number of these multichains is 
equal to the right-hand side of 
Equation~(\ref{eq:SP}) and the proof follows. 
\end{proof} 

We now present the final main result of this paper. 
We recall that the $h^\ast$-polynomial of an 
$n$-dimensional polytope in $\RR^n$ with vertices 
in $\ZZ^n$ is defined by the equation
\begin{equation} \label{eq:hstar-def} 
\sum_{m \ge 0} \Ehr(\qQ, m) t^m = 
\frac{h^\ast(\qQ, t)}{(1-t)^{n+1}}. 
\end{equation}

Given a word $w = (w_1, w_2,\dots,w_n) \in [n]^n$, 
we denote by $\asc(w)$ the number of indices 
(called \emph{ascents}) $i \in [n]$ such that 
$w_{i-1} < w_i$, where $w_0 := 1$.
\begin{proposition} \label{prop:h-star} 
We have
\begin{equation} \label{eq:h-star}
h^\ast(\qQ_\sigma, t) = \sum_{w \in \wW_\sigma} 
t^{\asc(w)},
\end{equation}
where $\wW_\sigma$ consists of all words $w \in 
[n]^n$ which have at least $s_i$ entries less than 
or equal to $s_i$ for every $i \in [k]$. 
Equivalently, the unique nondecreasing 
rearrangement $(u_1, u_2,\dots,u_n)$ of $w$ should 
satisfy $u_{s_i} \le s_i$ for every $i \in [k]$.
\end{proposition}

The proof employs the basics of the theory of 
EL-shellability of posets~\cite{Bj80}
\cite[Lecture~3]{Wa07}, which we briefly recall.
Let $(P, \preceq)$ be a finite graded poset 
which has a minimum element 
$\hat{0}$ and a maximum element $\hat{1}$ and
let $\eE(P)$ be the set of covering relations 
of $P$. An \emph{edge labeling} of $P$ (with 
integer labels) is any map $\lambda: \eE(P) 
\to \ZZ$. Let $C$ be an unrefinable chain $p_0 
\prec p_1 \prec \cdots \prec p_m$ of elements 
of $P$, so that $(p_{i-1}, p_i) \in \eE(P)$ for 
every $i \in [m]$. We define the \emph{label} 
of $C$ with respect to $\lambda$ as the sequence
\[ \lambda(C) = (\lambda(p_0, p_1), 
   \lambda(p_1, p_2),\dots,\lambda(p_{m-1}, 
	         p_m)) \]
and the \emph{descent set} of $C$ as 
$\Des(C) = \{ i \in [m-1]: \lambda(p_{i-1}, 
p_i) > \lambda(p_i, p_{i+1}) \}$. We say that 
$C$ is \emph{rising} with respect to $\lambda$ 
if $\lambda(C)$ is weakly increasing and that 
$C$ is \emph{lexicographically smaller} (with 
respect to $\lambda$) than an unrefinable chain 
$C'$ in $P$ of the same length if $\lambda(C)$ 
precedes $\lambda(C')$ in the lexicographic 
order. An edge labeling $\lambda$ of $P$ is 
called an \emph{EL-labeling} if for every 
nonsingleton interval $[u, v]$ in $P$: (i) 
there is a unique rising (with respect to 
$\lambda$) maximal chain in $[u, v]$; and (ii) 
this chain is lexicographically smallest among 
all maximal chains in $[u, v]$. Then 
\begin{equation} \label{eq:h-order}
h(\Delta(P),t) = \sum_{C \in \mM(P)} t^{\des(C)}, 
\end{equation}
where $\mM(P)$ is the set of all maximal 
chains of $P$ and $\des(C) = \# \Des(C)$ is the
number of descents of $C$ for $C \in \mM(P)$. 
Here $\Delta(P)$ stands for the \emph{order 
complex} of $P$ (meaning, the simplicial 
complex consisting of all chains in $P$) and the 
$h$-polynomial of $\Delta(P)$ is defined as 
$h(\Delta(P),t) = \sum_{F \in \Delta(P)} 
t^{\# F} (1-t)^{d+1 - (\# F)}$, where $d$ is the
rank of $P$.

\begin{proof}
[Proof of Proposition~\ref{prop:h-star}]
Let us write $\tau = \rev(\sigma)$. By 
Equation~(\ref{eq:hstar-def}) applied to 
$\qQ_\sigma$ and since 
\[ \sum_{m \ge 0} \zZ(P_\tau, m+1) t^m = 
   \frac{h(\Delta(P_\tau), t)} {(1-t)^{n+1}} \]
(see, for instance, \cite[Equation~(8)]{ADK24}) 
and $\Ehr(\qQ_\sigma, t) = \zZ(P_\tau, t+1)$ by 
Theorem~\ref{thm:EZ}, we have $h^\ast(\qQ_\sigma, 
t) = h(\Delta(P_\tau), t)$. Thus, it suffices to 
show that $h(\Delta(P_\tau), t)$ is equal to the 
sum on the right-hand side of 
Equation~(\ref{eq:h-star}).

For that reason, we consider the following 
edge labeling of the poset $\widehat{P}_\tau = 
P_\tau \cup \{ \hat{1} \}$, obtained from $P_\tau$
by adding a maximum element $\hat{1}$. We recall 
that $P_\tau$ is the subposet of the product order 
$\NN^n$ consisting of all points $(a_1, 
a_2,\dots,a_n) \in \NN^n$ such that
\[ a_{s_i + 1} + a_{s_i + 2} + \cdots + a_n \le 
   r_{i+1} + r_{i+2} + \cdots + r_k = n - s_i 
	 \]
for $0 \le i < n$. The cover 
relations are the pairs $(u,v) \in P_\tau \times
P_\tau$ such that $v$ can be obtained from $u$ 
by adding 1 in one of its coordinates. Given 
such a pair we set $\lambda(u,v) = i$, if $u$ 
and $v$ differ in the $i$th coordinate. We also 
set $\lambda(u, \hat{1}) = 1$ for every maximal
element $u \in P_\tau$ and claim that the 
resulting edge labeling $\lambda: 
\eE(\widehat{P}_\tau) \to \ZZ$ is an EL-labeling. 
Given a nonsingleton closed interval $[u, v]$ of 
$\widehat{P}_\tau$, we need to verify
that there is a unique rising (with 
respect to $\lambda$) maximal chain in $[u, v]$ 
and that this chain is lexicographically least 
among all maximal chains of $[u, v]$. This holds
in case $v \ne \hat{1}$ because then the labels
of the maximal chains of $[u, v]$ are exactly 
the permutations of a multiset; it holds if $v = 
\hat{1}$ as well since then there is a unique 
maximal chain with label $(1, 1,\dots,1)$ and 
this chain is the only rising maximal chain of 
$[u, v]$.

Since $\lambda$ is an EL-labeling, by 
Equation~(\ref{eq:h-order}) we have 
\[ h(\Delta(P_\tau), t) = \sum_w t^{\des(w)}, \]
where $w = (w_n,\dots,w_1, w_0)$ ranges over the 
labels of all maximal chains of 
$\widehat{P}_\tau$ and $\des(w)$ is the number of
indices $i \in [n]$ such that $w_i > w_{i-1}$. We 
observe that $w$ is such a label if and only if 
$w_0 = 1$ and at least $s_i$ entries of 
$(w_n,\dots,w_2, w_1)$ are less than or equal to 
$s_i$ for every $i \in [k]$ and the proof follows.
\end{proof}

\begin{remark} \label{rem:xue} \rm
Proposition~\ref{prop:h-star} specializes to a 
well known combinatorial interpretation of the 
$h^\ast$-polynomial (equivalently, of the 
$h$-polynomial of the $n$-fold edgewise 
subdivision) of the standard $n$-dimensional 
simplex when $\sigma = (n)$ has one part; see, 
for instance, \cite[Corollary~8]{SS12} for a 
generalization and refinement.
A generalization of Proposition~\ref{prop:h-star},
valid for the polytopes in $\RR^n$ defined by the 
equations 
\[ x_1 + x_2 + \cdots + x_{s_i} \le c_1 + c_2 + 
   \cdots + c_i \]

\medskip
\noindent
for $1 \le i \le k$, where $c_1, c_2,\dots,c_k$ are 
positive integers, has been found by Xue 
Yan~\cite{Yan25}.
\end{remark}

\medskip
\noindent
\textbf{Acknowledgements}. The author wishes to 
thank Frederic Chapoton for his crucial input to 
the proof of Conjecture~\ref{conj:poly}, described 
in Remark~\ref{rem:Cha}, and an anonymous reviewer 
for carefully reading the manuscript.



\end{document}